\documentclass{amsart}
\usepackage{graphicx}
\usepackage{amssymb}
\usepackage{amsfonts}
\usepackage{epsfig}
\swapnumbers
\newtheorem{theorem}{Theorem}[section]

\newtheorem{corollary}[theorem]{Corollary}

\theoremstyle{definition}

\numberwithin{equation}{section}
 \theoremstyle{plain}    
 
 \numberwithin{equation}{section} 
 \numberwithin{figure}{section} 
 \theoremstyle{plain}    
 \theoremstyle{plain}    
 \theoremstyle{remark}    
 \newtheorem*{acknowledgement*}{Acknowledgement} 

\newcommand{\cF}{{\mathcal F}}
\newcommand{\cG}{{\mathcal G}}
\newcommand{\cH}{{\mathcal H}}

\newcommand{\cP}{{\mathcal P}}

\newcommand{\cW}{{\mathcal W}}

\newcommand{\te}{{\theta}}

\newcommand{\vt}{{\vartheta}}

\newcommand{\vsig}{{\varsigma}}
\newcommand{\Om}{{\Omega}}
\newcommand{\om}{{\omega}}
\newcommand{\ve}{{\varepsilon}}
\newcommand{\del}{{\delta}}

\newcommand{\gam}{{\gamma}}
\newcommand{\Gam}{{\Gamma}}
\newcommand{\vf}{{\varphi}}

\newcommand{\up}{{\upsilon}}

\newcommand{\sig}{{\sigma}}
\newcommand{\al}{{\alpha}}
\newcommand{\be}{{\beta}}
\newcommand{\ka}{{\kappa}}
\newcommand{\la}{{\lambda}}

\newcommand{\bbN}{{\mathbb N}}

\newcommand{\bbR}{{\mathbb R}}

\newcommand{\bbZ}{{\mathbb Z}}

\newcommand{\frb}{{\mathfrak b}}

\begin{document}
\title[]{Erd\H os-R\' enyi Law of Large Numbers in the Averaging Setup.}%
 \vskip 0.1cm 
 \author{ Yuri Kifer\\
\vskip 0.1cm
Institute of Mathematics\\
Hebrew University\\
Jerusalem, Israel}%
\address{
 Institute of Mathematics, The Hebrew University, Jerusalem 91904, Israel}%
\email{ kifer@math.huji.ac.il}%

\thanks{}
\subjclass[2000]{Primary: 60F15 Secondary: 60F10, 34K33, 37D20, 37D35, 60J99}%
\keywords{laws of large numbers, large deviations, averaging, hyperbolic 
diffeomorphisms and flows, Markov processes.}%

\date{\today}
\begin{abstract}\noindent
We extend the Erd\H os--R\' enyi law of large numbers to the averaging
setup both in discrete and continuous time cases. We consider both stochastic 
processes and dynamical systems as fast motions whenever they are fast
mixing and satisfy large deviations estimates. In the continuous time case
we consider flows with large deviations estimates which allow a suspension
 representation and it turns out
that fast mixing of corresponding base transformations suffices for our results.
\end{abstract}
\maketitle
\markboth{Y.Kifer}{Averaging} 
\renewcommand{\theequation}{\arabic{section}.\arabic{equation}}
\pagenumbering{arabic}

\section{Introduction}\label{sec1}\setcounter{equation}{0}

Let $\xi_1,\xi_2,...$ be a sequence of independent identically distributed
(i.i.d.) random variables such that $E\xi_1=0$ and the moment generating
function $M(t)=Ee^{t\xi_1}$ exists. Denote by $I$ the Legendre transform 
of $\ln M$ and set $X_n=\sum_{m=1}^n\xi_m$ for $n\geq 1$ and $X_0=0$. 
The Erd\" os-R\' enyi law of large
 numbers from \cite{ER} says that with probability one,
\begin{equation}\label{1.1}
I(\be)\lim_{n\to\infty}\max_{0\leq m\leq n-[\frac {\ln n}{I(\be)}]}\frac
{X_{m+[\frac {\ln n}{I(\be)}]}-X_m}{\ln n}=\be
\end{equation}
for all $\be>0$ in some neighborhood of zero (actually, whenever 
$I(\be)<\infty$). A related result was proved earlier in \cite{Sh} where 
the maximum in (\ref{1.1}) is taken in $m$ up to $n$ and $\ln n$ is replaced 
in the above numerator by $\ln m$ so that this type of limits are sometimes 
called Erd\H os-R\' enyi-Shepp laws. There were numerous specifications and 
extensions  of the Erd\H os-R\' enyi law for the last 40 years (see, for
instance \cite{DD} and references there) while in \cite{No} a corresponding
limit law was derived. As the 
original version (\ref{1.1}) most of these results were valid only in the one 
 dimensional case, i.e. for random variables and not for random vectors. On
 the other hand, a functional form of (\ref{1.1}) suggested in \cite{Bor}
  holds true for i.i.d. random vectors, as well. 
 More recently papers on the Erd\H os-R\' enyi law appeared in the dynamical
 systems framework where extensions from i.i.d. to weakly dependent summands
 became necessary. In \cite{CC} the Erd\H os-R\' enyi law was derived for
 functions of iterates of expanding maps of the interval while an 
 extension of (\ref{1.1}) to stationary $\al$-mixing sequences was obtained 
 in \cite{DK} and to functions of some nonuniformly expanding dynamical 
 systems in \cite{DN}. In somewhat different direction an Erd\H os-R\' enyi
 law for Gibbs measure was derived earlier in \cite{Co}.
 
 In this paper we extend the Erd\H os-R\' enyi type results to the averaging
 setup which was not considered before generalizing all previous approaches to
 the problem (except for the case of nonconventional sums studied in 
 \cite{Ki4}). We consider the slow motion $X^\ve$ in both the discrete time case 
 \begin{equation}\label{1.2}
 X^\ve_{n+1}=X^\ve_n+\ve B(X^\ve_n,\xi_n),\quad X_0^\ve=x
 \end{equation}
 and in the continuous time case
 \begin{equation}\label{1.3}
 \frac {dX^\ve_t}{dt}=\ve B(X^\ve_t,\xi_t),\quad X_0^\ve=x.
 \end{equation}
 Here the fast motion $\xi_n,\, n\in\bbZ$ or $\xi_t,\, t\in\bbR$ is a stationary 
 stochastic process, in particular, it can be generated by a dynamical system
 $\xi_n=\xi_n(x)=f^nx,\, n\in\bbZ$ or $\xi_t=\xi_t(x)=f^tx,\ t\in\bbR$ 
 preserving some
 probability measure $\mu$ which plays the role of probability on the
 corresponding space where $x$ lives. In the discrete time we assume that
 $\xi_n$ is exponentially fast $\al$-mixing while in the continuous time
 case, in order to enable applications to important classes of dynamical
 systems, we asume that $\xi_t$ can be represented via so called suspension
 construction over an exponentially fast $\al$-mixing discrete time
 stationary process. 
 
 We observe that (\ref{1.2}) and (\ref{1.3}) are generalizations of usual 
 Ces\' aro averages of sums or integrals since if $B$ does not depend on 
 the slow motion $X^\ve$ then in the discrete time case
 \[
 X^\ve_{[1/\ve]}=x+\ve\sum_{0\leq n<[1/\ve]}B(\xi_n)\quad\mbox{and}\quad
 X^\ve_{1/\ve}=x+\ve\int_0^{1/\ve}B(\xi_t)dt
 \]
 in the continuous time case.
 
 If we fix an ergodic stationary measure $\mu$ of the process $\{\xi_n,\, 
 n\in\bbZ_+\}$ or $\{\xi_t,\, t\in\bbR_+\}$ then by the ergodic theorem the
 limits
 \begin{equation}\label{1.4}
 \lim_{N\to\infty}\frac 1N\sum_{n=0}^{N-1}B(x,\xi_n)=\bar B(x)\quad\mbox{or}
 \quad\lim_{T\to\infty}\frac 1T\int_0^TB(x,\xi_t)dt=\bar B(x)
 \end{equation}
 exist $\mu$-almost surely (a.s.) but, of course, they depend on $\mu$. The
 averaging principle proved rigorously about 70 years ago (see \cite{SVM} 
 and references there) says that if $B$ is Lipschitz continuous in the first
 variable then (\ref{1.4}) implies that
 \begin{equation}\label{1.5}
 \lim_{\ve\to 0}\sup_{0\leq t\leq T/\ve}|X^\ve_t-\bar X^\ve_t|=0
 \end{equation}
 where $\bar X^\ve_t$ is the averaged motion solving the equation
 \begin{equation}\label{1.6}
 \frac {d\bar X^\ve_t}{dt}=\ve\bar B(\bar X^\ve_t),\,\,\,\bar X_0^\ve=x.
 \end{equation}
 
 Since, in view of the above, the averaging principle can be considered as
  a generalization of the ergodic theorem, which in the probabilistic language
  is, essentially, a law of large numbers, it would be natural to ask whether 
  the Erd\H os-R\' enyi law can be generalized to the averaging setup, as well.
  It should be clear from the beginning that the latter cannot be obtained in
  so general circumstances as the averaging principle itself since the
  Erd\H os-R\' enyi law strongly relies on large deviations which can be
  proved only for certain classes of processes.
  
  In this paper we derive the Erd\H os-R\' enyi law type results in a 
  functional form for the slow motion $X^\ve_t$ (both in discrete and
  continuous time cases) in place of sums of random variables in (\ref{1.1}).
  When $B$ in (\ref{1.2}) or in (\ref{1.3}) does not depend on the first 
  variable our results yield the functional form of the Erd\H os-R\' enyi
  law which was introduced in \cite{Bor} in the particular case of sums
  of i.i.d. random vectors. When, in addition, $B$ is one dimensional this
  implies the Erd\H os-R\' enyi law in the form (\ref{1.1}) but in a much
  more general situation.
  
  This paper extends previous results on the Erd\H os-R\' enyi law in several 
  directions. First, it was never considered before in the averaging setup.
  Secondly, its functional form appeared before only for sums of i.i.d. vectors.
  Thirdly, this law was never dealt with in the continuous time case. Finally,
  we require weaker $\al$-mixing and not $\psi$-mixing conditions which
  appeared in \cite{DK}.

  Our results are applicable to several types of stationary processes $\xi_t$.
  On the probabilistic side $\xi_t$ can be, in particular, a Markov chain 
  satisfying an appropriate Doeblin condition or a nondegenerate diffusion
  process on a compact manifold. On the dynamical systems side we can take
  $\xi_t(x)=f^tx$ with $f$ being an Axiom A diffeomorphism or flow on
  a hyperbolic set (and could be considered also in a neghborhood of an
  attractor) while in the discrete time case additional options are 
  possible such as mixing subshifts of finite type, expanding transformations 
  and some maps of the interval.
  
  The structure of this paper is the following. In the next section we present
  precisely our general setup and formulate main results under conditions
  which include both probabilistic and dynamical systems examples mentioned
  above. In Sections \ref{sec3} and \ref{sec4} we give proofs of our results
  in the discrete and continuous time cases, respectively. In Appendix we
  discuss applications to various specific classes of stochastic processes
  and dynamical systems and recall properties of rate functions of large
  deviations needed in the proofs of our results.

\section{Preliminaries and main results}\label{sec2}\setcounter{equation}{0}

Let $M$ be a Polish (complete separable metric) space, $\bbR^d$ be a 
$d$-dimensional Euclidean
space and a bounded Borel map $B:\,\bbR^d\times M\to\bbR^d$ satisfies
\begin{equation}\label{2.1}
|B(x,y)-B(z,y)|\leq L_1|x-z|,\,\, |B(x,y)|\leq L_1
\end{equation}
for some $L_1>0$, all $x,z\in\bbR^d$ and any $y\in M$. We consider also 
a stationary ergodic stochastic process $\xi_t$ with discrete $t\in\bbZ_+=
\{ 0,1,...\}$ or continuous $t\in\bbR_+=\{ s\geq 0\}$ time on a probability
 space $(\Om,\cF,P)$ with values in $M$. Our setup includes also a sequence
 $\cF_{m,n}\subset\cF,\,-\infty\le m\le n\leq\infty$ of $\sig$-algebras such
 that $\cF_{m,n}\subset\cF_{m_1,n_1}$ whenever $m_1\le m$ and $n_1\ge n$ which 
 satisfies an exponentially fast $\al$-mixing condition (see, for instance,
 \cite{Bra}),
 \begin{equation}\label{2.2}
 \al(n)=\sup\big\{\big\vert P(A\cap B)-P(A)P(B)\big\vert :\, A\in
 \cF_{-\infty,k},\, B\in\cF_{k+n,\infty}\big\}\leq\ka_1^{-1}e^{-\ka_1n}
 \end{equation}
 for some $\ka_1>0$ and all $k,n\geq 0$.
 
 In the discrete time case we rely on the following approximation condition
 \begin{equation}\label{2.3}
 \zeta(n)=E\sup_{x,k}|B(x,\xi_k)-E(B(x,\xi_k)|\cF_{k-n,k+n})|\leq\ka_2^{-1}
 e^{-\ka_2n}
 \end{equation}
 for some $\ka_2>0$ and all $n\geq 0$.
 
 Next, set $\bar B(x)=EB(x,\xi_0)$, $\bar Z_t=\bar X^\ve_{t/\ve}$, $B_t(y)
 =B(\bar Z_t,y)$, $\bar B_t=EB_t(\xi_0)$ and $G_t(y)=B_t(y)-\bar B_t$ where
 $\bar X^\ve_t$ satisfies (\ref{1.6}), and so
 \begin{equation}\label{2.4}
 \frac {d\bar Z_t}{dt}=\bar B(\bar Z_t).
 \end{equation}
 Define
 \begin{equation}\label{2.5}
 Y_{t,r}(u)=\frac 1r\sum_{0\leq j\leq ru}G_t(\xi_j),\,\, u\in[0,1],\, r\in\bbN
 =\{ 1,2,...\}.
 \end{equation}
 Denote the space of continuous curves $\gam:\,[0,1]\to\bbR^d$ by $C([0,1],
 \bbR^d)$ and assume that for any $\gam\in C([0,1],\bbR^d)$ and $t\in[0,T]$
 the limit 
 \begin{equation}\label{2.6}
 \lim_{r\to\infty}\frac 1r\ln E\exp(r\int_0^1(\gam_u,G_t(\xi_{[ru]}))du)=
 \int_0^1 \Pi_t(\gam_u)du
 \end{equation}
 exists, where $\Pi_t(\frb),\,\frb\in\bbR^d$ is a convex twice differentiable
 function such that $\nabla_{\frb}\Pi_t(\frb)|_{\frb=0}=0$ and the Hessian 
 matrix $\nabla^2_{\frb}\Pi_t(\frb)|_{{\frb}=0}$ is positively definite.
 Here $(\cdot,\cdot)$ denotes the inner product. Let
 \begin{equation}\label{2.7}
 I_t(\be)=\sup_{\frb}((\frb,\be)-\Pi_t(\frb))
 \end{equation}
 and for any $\gam\in C([0,1],\bbR^d)$ set
 \begin{equation}\label{2.8}
 S_t(\gam)=\int_0^1I_t(\dot\gam_s)ds,\quad\dot\gam_s=\frac {d\gam_s}{ds}
 \end{equation}
 if $\gam$ is absolutely continuous and $S_t(\gam)=\infty$ for otherwise.
 It follows from the above (see, for instance, Section 7.4 in \cite{FW})
 that $Y_{t,r}$ satisfies large deviations estimates in the form that
 for any $a,\del,\la>0$ and every $\gam\in C([0,1],\bbR^d)$, $\gam_0=0$ 
 there exists $r_0>0$ such that for $r\geq r_0$,
 \begin{eqnarray}\label{2.9}
 &P\{\rho(Y_{t,r},\gam)<\del\}\geq\exp(-r(S_t(\gam)+\la))\quad\mbox{and}\\
 &P\{\rho(Y_{t,r},\Phi_t(a))\geq\del\}\leq\exp(-r(a-\la))\nonumber
 \end{eqnarray}
 where $\rho(\gam,\eta)=\sup_{s\in[0,1]}|\gam_s-\eta_s|$ and
 $\Phi_t(a)=\{\gam\in C([0,1],\bbR^d):\,\gam_0=0,\, S_t(\gam)\leq a\}$.

 Since $S_t$ is a lower semi-continuous functional then each $\Phi_t(a)$ is
 a closed set and, moreover, it is compact for any finite $a$. Indeed, 
 $|\Pi_t(\frb)|\leq 2L_1|\frb|$ by (\ref{2.1}) and (\ref{2.6}) which implies
 by (\ref{2.7}) that $I_t(\be)=\infty$ provided $|\be|>2L_1$ (take 
 $\frb=a\be/|\be|$ in (\ref{2.7}) and let $a\to\infty$). Hence, $|\dot\gam_s|
 \leq 2L_1$ for Lebesgue almost all $s\in[0,1]$ if $\gam\in\Phi_t(a)$, and
 so the latter set is bounded and equicontinuous which by the Arzel\` a-Ascoli
 theorem implies its compactness.

 In the discrete time case (\ref{1.2}) for any $\ve>0,\, t\in[0,T],\,
 u\in[0,1]$ and $N\in\bbN$ set
 \begin{equation*}
 V^{\ve,N}_t(u)=\frac {X^\ve_{[t/\ve]+[b_t(\ve,N,u)]}-X^\ve_{[t/\ve]}}
  {\ve b_t(\ve,N)}-u\bar B_t
  \end{equation*}
  where 
  \[
  b_t(\ve,N,u)=c_{\tau(t,N)}u\ln\frac 1\ve,\, b_t(\ve,N)=b_t(\ve,N,1),\,
  \tau(t,N)=[Nt/T]T/N
  \]
  and $c_t$ is a function on $[0,T]$ such that $0<\hat c^{-1}\leq c_t\leq
  \hat c<\infty$ for some $\hat c$.

 \begin{theorem}\label{thm2.1}
 Assume that the conditions (\ref{2.2}), (\ref{2.3}) and (\ref{2.6}) hold
 true. Then $V^{\ve,N}_t(u)$ defined above satisfies
 \begin{equation}\label{2.10}
 \lim_{N\to\infty}\limsup_{\ve\to 0}sup_{0\leq t\leq T}\,
 \rho\big(V_t^{\ve,N},\Phi_{\tau(t,N)}(c^{-1}_{\tau(t,N)})\big)=0\quad 
 \mbox{a.s.}
 \end{equation}
 and 
 \begin{equation}\label{2.11}
 \lim_{N\to\infty}\limsup_{\ve\to 0}sup_{\gam\in\Phi_0(1/c_0)}
 \inf_{0\leq t\leq T}\rho(V_t^{\ve,N},\gam)=0 \quad\mbox{a.s.}
 \end{equation}
 \end{theorem}

 \begin{corollary}\label{cor2.2} Suppose that $B(x,y)=B(y)$ does not depend
 on the first variable. Then $G_t,\, Y_{t,r},\,\Pi_t,\, I_t$ and $S_t$ in 
 (\ref{2.5})--(\ref{2.8}) do not depend on $t$, and so $\Phi_t(a)=\Phi(a)$
 does not depend on $t$, as well. Let $c_t\equiv c>0$ be a constant then
 $V^{\ve,N}_t=V^\ve_t$ does not depend on $N$. Set $\cW^\ve_c=\cup_{0\leq t
 \leq T}V^\ve_t$. Then for any $c>0$,
 \begin{equation}\label{2.12}
 \lim_{\ve\to 0}H(\cW^\ve_c,\,\Phi(1/c))=0\quad\mbox{a.s.}
 \end{equation}
 where $H(\Gam_1,\Gam_2)=\inf\{\del>0:\,\Gam_1\subset \Gam^\del_2,\,
 \Gam_2\subset \Gam^\del_1\}$ is the Hausdorff distance between sets of
 curves with respect to the uniform metric $\rho$ (and $\Gam^\del=
 \{\gam:\,\rho(\gam,\Gam)<\del\}$ is the $\del$-neighborhood of $\Gam$).
 \end{corollary}
 
 \begin{corollary}\label{cor2.3} Suppose that $d=1$ and set 
 $c_t=\frac 1{I_t(\be)}$ where 
 $\be_t>\be>0$ and $\be_t=\sup\{\be>0:\, I_t(\be)<\infty\}$.
 Then 
 \begin{equation}\label{2.13}
 \lim_{N\to\infty}\lim_{\ve\to 0}\sup_{0<t<T}V^{\ve,N}_t(1)=\be\quad
 \mbox{a.s.}
 \end{equation}
 where $\lim_{N\to\infty}\lim_{\ve\to 0}=\lim_{N\to\infty}\limsup_{\ve\to 0}
 =\lim_{N\to\infty}\liminf_{\ve\to 0}$.
 In particular, if $d=1$ and $B(x,y)=B(y)$ does not depend on the first
 variable then $I_t\equiv I$ does not depend on $t$ and (\ref{2.13}) holds
 true with $c_t\equiv \frac 1{I(\be)}$ for all $0<\be<\be_0=
 \sup\{\be>0:\, I(\be)<\infty\}$.
 \end{corollary}
 
 We observe that the proof of Theorem \ref{thm2.1} and of Corollary
 \ref{cor2.3} require certain time discretization which cannot be achieved
  relying on some continuity properties since $I_t$ and $S_t$ are only lower
 semi continuous in $t$. By this reason we had to introduce $\tau(t,N)$ and 
 to have the second limit as $N\to\infty$. In fact, 
 if $|\be|$ is small enough then $I_t(\be)$ is continuous in $t$ but in order 
 to use this we would have to consider only curves $\gam$ with uniformly small
 speeds $|\dot\gam_t|$ which would not be a natural restriction.
 Though we work in a substantially more general averaging setup the strategy 
 of our proof of Theorem \ref{thm2.1} resembles previous works, in particular,
  \cite{DK} and \cite{DN} but observe that we rely only on $\al$-mixing and
  not on a stronger $\psi$-mixing assumed in the above papers. Large deviations
  estimates for hyper-geometrically fast $\al$ and $\phi$-mixing stationary 
  sequences were derived in \cite{Bry} and \cite{BD} while existence of such 
  processes follows 
  from Theorem 2 in \cite{Bra0}. We observe also that Corollary \ref{cor2.3}
  does not require full strength of large deviations in the form (\ref{2.9})
  and it suffices to have here usual level one large deviations estimates
  for $Y_{t,r}(1)$ in the form 
  \begin{eqnarray}\label{2.14}
  &\limsup_{r\to\infty}\frac 1r\ln P\{ Y_{t,r}(1)\in K\}\leq-\inf_{\frb\in K}
  I(\frb)\quad\mbox{and}\\
  &\liminf_{r\to\infty}\frac 1r\ln P\{ Y_{t,r}(1)\in U\}\geq-\inf_{\frb\in U}
  I(\frb)\nonumber
  \end{eqnarray}
  for any closed $K$ and open $U$ subsets of real numbers. Our method will
  still go through with minor modifications if the exponentially fast
  decay in (\ref{2.2}) and (\ref{2.3}) is replaced by a streched exponential
  one, i.e. by $exp(-\ka n^{\del})$ for some $\ka,\del>0$.
 
 Next, we deal with the continuous time case.  In addition to a stationary
 ergodic process $\xi_t,\, t\in\bbR_+$ on a probability space $(\Om,\cF,P)$
 with a path shift operator $\vt:\Om\to\Om$ we consider now an embedded discrete
 time process $\eta_k,\, k\in\bbZ_+$ related to $\xi_t(\om)=\xi_0(\vt^t\om)$
 by means of measurable  maps $\vf:\,\Om\to\hat\Om\subset\Om$, 
 $\te:\,\hat\Om\to\hat\Om$ and a  measurable function $\vsig:\,\hat\Om\to\bbR_+$
  such that $\vf^{-1}(\hat\om)\subset\{\vt^t\hat\om:\, 0\leq t<\vsig(\hat\om)\}$
  and for any $\hat\om=\vf(\om)$ and $k\geq 0$,
 \begin{equation}\label{2.15}
 \vt^{\vsig(\hat\om)}\hat\om=\te\hat\om,\,\eta_k(\hat\om)=\eta_0(\te^k\hat\om)
 \,\,\,\mbox{and}\,\,\,\eta_0(\hat\om)=\xi_0(\hat\om).
 \end{equation}
 Set $\hat P=\vf P$ and $\hat\cF=\vf\cF$. We assume that there exists a 
 $\te$-invariant probability measure $Q$ on $(\hat\Om,\hat\cF)$ equivalent
 to $\hat P$ and such that
 \begin{equation}\label{2.16}
 L^{-1}_2\leq\frac {dQ}{d\hat P}\leq L_2
 \end{equation}
 for some $L_2>0$. Thus, $\eta_k,\, k\geq 0$ is a stationary process with 
 respect to $Q$.
 For each $\hat\om=\vf(\om)$ and $j\geq 0$ set
 \[
 \hat B(x,\om)=\hat B(x,\hat\om)=\int_{0}^{\vsig(\hat\om)}B(x,\xi_u(\hat\om))du.
 \]
 We assume that there exist sub $\sig$-algebras $\cF_{m,n}\subset\hat\cF$ on 
 $\hat\Om$ satisfying (\ref{2.2}) with $Q$ in place of $P$ and such that for 
 some $\ka_3,L_3>0$ and all $n\geq 0$,
 \begin{equation}\label{2.17}
 \sup_{x,j}\hat E\big\vert\hat B(x,\cdot)\circ\te^j-\hat E(\hat B(x,\cdot)
 \circ\te^j|\cF_{j-n,j+n})\big\vert\leq\ka_3^{-1}e^{-\ka_3n},
 \end{equation}
 \begin{equation}\label{2.18}
 L_3^{-1}\leq\vsig\leq L_3\quad\mbox{and}\quad\sup_{j\geq 0}\hat E|\vsig\circ
 \te^j-\hat E(\vsig\circ\te^j|\cF_{j-n,j+n})|\leq\ka_3^{-1}e^{-\ka_3n}
 \end{equation}
 where $\hat E$ is the expectation with respect to $\hat P$. Our proof
 will also go through if in the approximation conditions (\ref{2.17}) and
 (\ref{2.18}) we take the expectation $E_Q$ with respect to $Q$ in place of 
 $\hat E$.
 We assume also an upper large deviations bound for $\sig_n=\sum_{i=0}^{n-1}
 \vsig\circ\te^i$ in the form
 that for any $\del>0$ there exists $\ka_\del>0$ such that
 \begin{equation}\label{2.19}
 Q\{\frac 1n\sig_n\geq\bar\vsig(1+\del)\}\leq\ka_\del^{-1}e^{-\ka_\del n}
 \end{equation}
 where $\bar\vsig=E_Q\vsig$.
 
 Now, let $X^\ve_t$ be defined by (\ref{1.3}), $\bar X^\ve_t$ be defined
 by (\ref{1.6}) and, again, $\bar Z_t=\bar X^\ve_{t/\ve},\,\bar B_t(y)=
 B(\bar Z_t,y),\,\bar B_t=EB_t(\xi_0)$ and $G_t(y)=B_t(y)-B_t$. Set
 \begin{equation}\label{2.20}
 Y_{t,r}(u)=\frac 1r\int_0^{ru}G_t(\xi_s)ds,\,\, u\in[0,1].
 \end{equation}
 Assume that 
 \begin{equation}\label{2.21}
 \lim_{r\to\infty}\frac 1r\ln E\exp(r\int_0^1(\gam_u,G_t(\xi_{ru}))du)=\int_0^1
 \Pi_t(\gam_u)du
 \end{equation}
 exists for any $\gam\in C([0,1],\bbR^d)$ with $Y_{t,r}$ defined by 
 (\ref{2.20}), where, again, $\Pi_t(\frb),\,\frb\in\bbR^d$ is a convex twice 
 differentiable
 function such that $\nabla_\frb\Pi_t(\frb)|_{\frb=0}=0$ and the Hessian matrix
 $\nabla^2_\frb\Pi_t(\frb)|_{\frb=0}$ is positively definite.
  Then, again, by \cite{FW} the large
 deviations estimates (\ref{2.9})  hold true with $S_t(\gam)$
 defined by (\ref{2.7}) and (\ref{2.8}).

 \begin{theorem}\label{thm2.4} Assume that the conditions (\ref{2.2}),
 (\ref{2.6}) and (\ref{2.16})--(\ref{2.19}) hold true. Then $V^\ve_t$
  defined by 
 \begin{equation*}
 V^{\ve,N}_t(u)=\frac {X^\ve_{t/\ve+b_t(\ve,N,u)}-X^\ve_{t/\ve}}
  {\ve b_t(\ve,N)}-u\bar B_t
  \end{equation*}
  satisfies (\ref{2.10}) and (\ref{2.11}).
 \end{theorem}

 \begin{corollary}\label{cor2.5}
 Corollaries \ref{cor2.2} and \ref{cor2.3} remain true under the conditions of 
 Theorem \ref{thm2.4}, as well.
 \end{corollary}
 
 The proof of the upper bound (\ref{2.10}) proceeds in the continuous time case
 similarly to Theorem \ref{thm2.1} since it uses essentially only stationarity
 of the process $\xi_t$ and the large deviations bounds which are our
 assumptions in the above setup. On the other hand, the proof of the lower
 bound (\ref{2.11}) requires additional ingradients in the continuous time
 case since in order to accomodate important classes of dynamical systems
 we do not impose strong mixing conditions on the process $\xi_t$ itself
 but only on the base discrete time process $\eta_k$ (via the family of
 $\sig$-algebras $\cF_{m,n}$).
 To the best of our knowledge the Erd\H os-R\' enyi law of large numbers type
 results were not obtained before in any continuous time framework. If 
 the fast motion $\xi_t$ is, say, a nondegenerate diffision process on a
 compact manifold then it is exponentially fast $\psi$-mixing with respect
 to $\sig$-algebras generated by itself (see, for instance \cite{Bra}) and
 in this case Theorem \ref{thm2.4} is easy to derive from Theorem \ref{thm2.1}.
 Our setup is adapted, in particular, to important classes of dynamical 
 systems such as Axiom A flows which may exhibit arbitrarily slow mixing
  but, on the other hand, can be represented by means of the above suspension
  construction over an exponentially fast $\psi$-mixing base transformation
   (see \cite{BR}).

 \section{Discrete time case}\label{sec3}\setcounter{equation}{0}
\subsection{Basic estimates}\label{subsec3.1}
Taking into account (\ref{1.2}) write 
\begin{eqnarray}\label{3.1}
&\ve^{-1}(X^\ve_{[t/\ve]+[b_t(\ve,N,u)]}-X^\ve_{[t/\ve]})=\sum_{[t/\ve]< n
\leq [t/\ve]+[b_t(\ve,N,u)]}B(X^\ve_n,\xi_n)\\
&=\sum_{[t/\ve]<n\leq [t/\ve]+[b_t(\ve,N,u)]}B_t(\xi_n)+\Psi^{(1)}_{\ve,N}
(t,u)+\Psi^{(2)}_{\ve,N}(t,u)\nonumber
\end{eqnarray}
where
\[
\Psi^{(1)}_{\ve,N}(t,u)=\sum_{[t/\ve]<n\leq [t/\ve]+[b_t(\ve,N,u)]}
(B(X^\ve_n,\xi_n)-B(\bar X^\ve_n,\xi_n))
\]
and
\[
\Psi^{(2)}_{\ve,N}(t,u)=\sum_{[t/\ve]<n\leq [t/\ve]+[b_t(\ve,N,u)]}
(B(\bar X^\ve_n,\xi_n)-B_t(\xi_n)).
\]
By (\ref{2.1}),
\begin{equation}\label{3.2}
|\Psi^{(1)}_{\ve,N}(t,u)|\leq L_1b_t(\ve,N)\Psi^{(3)}_{\ve,N}
\end{equation}
where by the averaging principle (\ref{1.5}) with probability one uniformly
 in $N$,
\begin{equation}\label{3.3}
\Psi_{\ve,N}^{(3)}=\sup_{0\leq n\leq T(\ve^{-1}-\ln\ve)}|X^\ve_n-\bar X^\ve_n|
\to 0\,\,\,\mbox{as}\,\,\,\ve\to 0.
\end{equation}
Furthermore, by (\ref{1.4}), (\ref{1.6}) and (\ref{2.1}),
\begin{equation}\label{3.4}
|\Psi^{(2)}_{\ve,N}(t,u)|\leq L_1b_t(\ve,N)\Psi^{(4)}_{\ve,N}
\end{equation}
where
\begin{equation}\label{3.5}
\Psi_{\ve,N}^{(4)}=\sup_{0\leq t\leq T}\,\sup_{[t/\ve]< n\leq[t/\ve]+
b_t(\ve,N)}|\bar X^\ve_n-\bar X^\ve_{t/\ve}|\leq L_1\ve b_t(\ve,N)
\to 0\,\,\,\mbox{as}\,\,\,\ve\to 0.
\end{equation} 

It follows from (\ref{3.1})--(\ref{3.5}) that (\ref{2.10}) and (\ref{2.11})
will hold true provided we establish these limits for
\[
\hat V_t^{\ve,N}(u)=\frac 1{b_t(\ve,N)}\sum_{[t/\ve]< n\leq[t/\ve]+
[b_t(\ve,N,u)]}(B_t(\xi_n)-\bar B_t),\,\, t\in[0,T],\, u\in[0,1]
\]
in place of $V_t^{\ve,N}(u)$. Next, observe that
\begin{equation}\label{3.6}
|B_t(y)-\bar B_{\tau(t,N)}|+|\bar B_t-\bar B_{\tau(t,N)}|\leq 2L_1|\bar Z_t-
\bar Z_{\tau(t,N)}|\leq 2L_1^2|t-\tau(t,N)|\leq 2L_1^2TN^{-1}.
\end{equation}
It follows from here that (\ref{2.10}) and (\ref{2.11}) will hold true 
provided we obtain these limits for
\begin{equation}\label{3.7}
W_t^{\ve,N}(u)=\frac 1{b_t(\ve,N)}\sum_{[t/\ve]< n\leq[t/\ve]+
[b_t(\ve,N,u)]}(B_{\tau(t,N)}(\xi_n)-\bar B_{\tau(t,N)}),
\end{equation}
$t\in[0,T],\, u\in[0,1]$ in place of $V_t^{\ve,N}(u)$.

\subsection{The upper bound}\label{subsec3.2}
For $k=0,...,N,\, l=0,1,...,[T/\ve]$ and $u\in[0,1]$ set
\[
\hat W^{\ve,N}_{k,l}(u)=b^{-1}_{kT/N}(\ve,N)\sum_{l\leq n\leq l+
[b_{kT/N}(\ve,N)]}(B_{kT/N}(\xi_n)-\bar B_{kT/N}).
\]
In order to prove (\ref{2.10}) for $W^{\ve,N}_t$ given by (\ref{3.7}) in
place of $V_t^{\ve,N}$ it suffices to show that with probability one,
\begin{equation}\label{3.8}
\limsup_{N\to\infty}\limsup_{\ve\to 0}\max_{0\leq k\leq N}
\max_{0\leq l\leq[T/\ve]}\rho(\hat W_{k,l}^{\ve,N},\,\Phi_{kT/N}(1/c_{kT/N})=0.
\end{equation}
For any $\ve,\del>0$ and each integer $N\geq 1$ define the event
\[
\Gam_\del(\ve,N)=\{\max_{0\leq k\leq N}\max_{0\leq l\leq[T/\ve]}\rho
(\hat W^{\ve,N}_{k,l},\,\Phi_{kT/N}(1/c_{kT/N}))>\del\}.
\]
Then
\begin{equation}\label{3.9}
P(\Gam_\del(\ve,N))\leq\sum_{0\leq k\leq N}\sum_{0\leq l\leq[T/\ve]}
P(\hat\Gam_{\del,k,l}(\ve,N))
\end{equation}
where 
\[
\hat\Gam_{\del,k,l}(\ve,N)=\{\rho(\hat W^{\ve,N}_{k,l},\,
\Phi_{kT/N}(1/c_{kT/N}))>\del\}.
\]

Recall, that each $\Phi_t(a)$ is a compact set.
It follows that for any $t\geq 0$ and $a,\del>0$ there exists $\sig=
\sig_{t,a,\del}$ such that (cf. \cite{Bor}),
\begin{equation}\label{3.10}
\Phi_t(a+\sig)\subset U_\del(\Phi_t(a))=\{ x:\,\rho(x,\,\Phi_t(a))<\del\}.
\end{equation}
Indeed, if for some $t\geq 0$ and $a,\del>0$ the sets $Q_\sig=\Phi_t(a+\sig)
\setminus U_\del(\Phi_t(a))\ne\emptyset$ for all $\sig>0$ then $\cap_{\sig>0}
Q_\sig\ne\emptyset$ since $Q_\sig,\,\sig>0$ are compact and $Q_\sig\supset
Q_{\sig'}$ when $\sig>\sig'$. But if $\gam_0\in\cap_{\sig>0}Q_\sig$ then
$\gam_0\in\Phi_t(a)$ which contradicts the fact that $\gam_0\not\in U_\del
(\Phi_t(a))$.

Observe that
\[
\hat W^{\ve,N}_{k,0}(u)=Y_{kTN^{-1},c_{kTN^{-1}}\ln(1/\ve)}(u)
\]
with $Y_{t,r}$ defined by (\ref{2.5}). Hence, choosing $\sig>0$ satisfying 
(\ref{3.10}) for $a=c^{-1}_{kT/N}$, $t=kT/N$ and $\del>0$ we obtain
employing the upper bound of large deviations from (\ref{2.9}) that
\begin{eqnarray}\label{3.11}
& P(\hat\Gam_{2\del,k,l}(\ve,N))\leq P\{\rho(\hat W^{\ve,N}_{k,l},\,
\Phi_{kT/N}(c^{-1}_{kT/N}+\sig))>\del\}\\
&\leq\exp(-b_{kT/N}(\ve)(c^{-1}_{kT/N}+\sig-\la))\leq\ve^{1+\hat\sig c_{kT/N}}
e^{d_{k,N}}\nonumber
\end{eqnarray}
where $\ve\leq\ve_0$ and $\ve_0>0$ is chosen so small that (\ref{3.11}) holds
true for some $\la<\sig$, and so $\hat\sig=\hat\sig_{k,\del}=\sig-\la>0$ while 
$d_{k,N}
=c^{-1}_{kT/N}+\hat\sig_{k,\del}$. Here $\hat\sig$ depends on $k$ and $\del$ so 
we take $\hat\sig_\del=\min_{0\leq k\leq N}\hat\sig_{k,\del}$. Now, by 
(\ref{3.9}) and (\ref{3.10}),
\begin{equation}\label{3.12}
P(\Gam_{2\del}(\ve,N))\leq T\sum_{0\leq k\leq N}\ve^{\hat\sig_{k,\del}c_{kT/N}}
e^{d_{k,n}}\leq TNe^{\hat c^{-1}+\hat\sig_\del}\ve^{\hat\sig_\del\hat c^{-1}}
\end{equation}
where, recall, $\hat c^{-1}\geq c_t\geq\hat c>0$.

Choose the sequence $\ve_n=n^{-\frac {2\hat c}{\hat\sig_\del}}$. Then 
(\ref{3.12}) together with the Borel-Cantelli lemma yield that for $P$-almost
all $\om$ there exists $n_\del(\om)<\infty$ such that $\om\not\in\Gam_{2\del}
(\ve_n,N)$ for all $n\geq n_\del(\om)$. It follows that with probability one,
\begin{equation}\label{3.13}
\limsup_{n\to\infty}\max_{0\leq k\leq N}\max_{0\leq l\leq[T/\ve_n]}
\rho(W^{\ve_n,N}_{k,l},\,\Phi_{kT/N}(c^{-1}_{kT/N}))\leq 2\del.
\end{equation}
Now observe that if $\ve_n\leq\ve<\ve_{n-1}$ then
\[
|b_t(\ve,N,u)-b_t(\ve_n,N,u)|\leq 2\hat c^2\hat\sig_\del^{-1}\ln\frac n{n-1}
\to 0\,\,\,\mbox{as}\,\,\, n\to\infty.
\]
This together with (\ref{3.13}) gives that with probability one,
\begin{equation*}
\limsup_{\ve\to 0}\max_{0\leq k\leq N}\max_{0\leq l\leq[T/\ve]}
\rho(W^{\ve,N}_{k,l},\,\Phi_{kT/N}(c^{-1}_{kT/N}))\leq 2\del
\end{equation*}
for any $N\in\bbN$. Since $\del>0$ is arbitrary we obtain (\ref{3.8}) even
without $\limsup_{N\to\infty}$ yielding (\ref{2.10}).

\subsection{The lower bound}\label{subsec3.3}
In view of (\ref{3.1})--(\ref{3.6}) it suffices to establish (\ref{2.11})
for $W^{\ve,N}_t$ in place of $V_t^{\ve,N}$. Next, we have
\begin{eqnarray}\label{3.14}
&P\big\{\inf_{0\leq t\leq T}\rho(W^{\ve,N}_t,\,\gam)\geq\del\}\leq 
P\{\cap_{0\leq t
<T/N}\{\rho(W^{\ve,N}_t,\,\gam)\geq\del\}\big\}\\
&=P\big\{\cap_{0\leq l<[T/N\ve]}\{\rho(\hat W^{\ve,N}_{0,l},\,\gam)\geq\del\}
\big\}
\nonumber\\
&\leq P\big\{\cap_{0\leq j<[T/N\ve][\ln^2\ve]^{-1}}\{\rho(\hat 
W^{\ve,N}_{0,j[\ln^2\ve]},\,\gam)\geq\del\}\big\}\overset{def}=I_\ve.\nonumber
\end{eqnarray}

Now, we are going to rely on mixing and approximation assumptions (\ref{2.2})
and (\ref{2.3}). Set
\[
\tilde W_l^{\ve,N}(u)=E(\hat W_{0,l}^{\ve,N}(u)\vert\cF_{l-[\frac 13\ln^2\ve],
l+[\frac 13\ln^2\ve]}),\, u\in[0,1].
\]
Then by (\ref{2.3}) for all $\ve>0$ small enough,
\begin{equation}\label{3.15}
\max_{0\leq l\leq[T/\ve]}E\sup_{u\in[0,1]}|\hat W_{0,l}^{\ve,N}(u)-
\tilde W_l^{\ve,N}(u)|\leq\exp(-\frac {\ka_2}4\ln^2\ve)
\end{equation}
where we use that for any random vector $\Xi$ and $\sig$-algebras 
$\cG\subset\cH$,
\[
E|\Xi-E(\Xi|\cH)|\leq E|\Xi-E(\Xi|\cG)|+E|E(\Xi|\cG)-E(\Xi|\cH)|\leq
2E|\Xi-E(\Xi|\cG)|.
\]

Set
\[
J_\ve=P\big\{\cap_{0\leq j<[T/N\ve][\ln^2\ve]^{-1}}\{\rho(\tilde 
W^{\ve,N}_{0,j[\ln^2\ve]},\,\gam)\geq\del/2\}\big\}.
\]
Then by (\ref{3.15}) and the Chebyshev inequality,
\begin{eqnarray}\label{3.16}
&I_\ve\leq J_\ve+\sum_{0\leq j<[T/N\ve][\ln^2\ve]^{-1}}P\big\{\rho(
\hat W^{\ve,N}_{0,j[\ln^2\ve]},\tilde W^{\ve,N}_{0,j[\ln^2\ve]})\geq\del/2\big\}
\\
&\leq J_\ve+2[T/N\ve][\ln^2\ve]^{-1}\del^{-1}\exp (-\frac {\ka_2}4\ln^2\ve).
\nonumber\end{eqnarray}

Next, we use (\ref{2.2}) which yields easily by iduction that for any events 
$A_1,...,A_k$ such that  $A_i\in\cF_{m_i,n_i}$, where $m_i\leq n_i<n_i+l_i
\leq m_{i+1}$ for all $i=1,...,k$ with $m_{k+1}=\infty$, 
 \begin{equation}\label{3.17}
 P(\cap_{1\leq i\leq k}A_i)\leq \prod_{1\leq i\leq k}P(A_i)+
 \sum_{1\leq i\leq k-1}\al(l_i).
 \end{equation}
 Indeed, (\ref{3.17}) follows for $k=2$ directly from (\ref{2.2}). If 
 (\ref{3.17}) holds true for $k-1$ in place of $k$ then applying (\ref{2.2})
 to $A=\cap_{1\leq i\leq k-1}A_i$ and $B=A_k$ we derive (\ref{3.17}).
 Now, taking into account that the random vectors 
 $\tilde W^{\ve,N}_{0,j[\ln^2\ve]}$ are $\cF_{j[\ln^2\ve]-[\frac 13\ln^2\ve],
 \, j[\ln^2\ve]+[\frac 13\ln^2\ve]}$-measurable we obtain from (\ref{2.2}) 
 and (\ref{3.17}) that for all $\ve>0$ small enough,
 \begin{equation}\label{3.18}
 J_\ve\leq\prod_{0\leq j\leq [T/N\ve][\ln^2\ve]^{-1}} P\{\rho(
 \tilde W^{\ve,N}_{0,j[\ln^2\ve]},\gam)\geq\del/2\}+
 \ve^{-1}\exp(-\frac {\ka_1}4\ln^2\ve).
 \end{equation}
 Using (\ref{3.15}) and the Chebyshev inequality again we have that
 \begin{eqnarray}\label{3.19}
 &P\{\rho(\tilde W^{\ve,N}_{l},\gam)\geq\del/2\}\leq P\{\rho(\hat 
 W^{\ve,N}_{0,l},\gam)\geq\del/4\}\\
 &+P\{\rho(\tilde W^{\ve,N}_{l},\hat W^{\ve,N}_{0,l})\geq\del/4\}\leq P\{\rho(
 \hat W^{\ve,N}_{0,l},\gam)\geq\del/4\}+\frac 4\del\exp(-\ka_2\ln^2\ve).
 \nonumber\end{eqnarray}
 
 Since $\gam\in\Phi_0^{1/c_0}$ then $I_0(\dot\gam_u)<\infty$ for Lebesgue 
 almost all $u\in[0,1]$ and, as explained in Appendix, if $I_0(\be)<\infty,\,
 \be\in\bbR^d$ then $I_0(a\be)<I_0(\be)$ for $0<a<1$. Hence, if we define
 \[
 \eta_u=(1-\del(8\sup_{v\in[0,1]}|\gam_u|)^{-1})\gam_u,\, u\in[0,1]
 \]
 then
 \begin{equation}\label{3.20}
 \rho(\gam,\eta)\leq\del/8\quad\mbox{and}\quad S_0(\eta)\leq S_0(\gam)-a\leq
 \frac 1{c_0}-a
 \end{equation}
 for some $a$.
 
 Next, we write
 \begin{equation}\label{3.21}
P\{\rho(\hat W^{\ve,N}_{0,l},\gam)\geq\del/4\}\leq P\{\rho(\hat 
W^{\ve,N}_{0,l},\eta)\geq\del/8\}=1-P\{\rho(\hat W^{\ve,N}_{0,l},\eta)<\del/8\}.
\end{equation}
Now, taking into account stationarity of the process $\xi_t$ and that 
$\hat W^{\ve,N}_{0,0}=Y_{0,c_0\ln\frac 1\ve}$
(with the latter defined by (\ref{2.5})) and relying on the lower large 
deviations bound in (\ref{2.9}) we obtain that for any $\del,\la>0$ there 
exists $\ve_0>0$ such that for all positive $\ve<\ve_0$,
 \begin{eqnarray}\label{3.22}
 &P\{\rho(\hat W^{\ve,N}_{0,l},\eta)<\frac \del 8\}=P\{\rho(\hat W^{\ve,N}_{0,0},
 \eta) <\frac \del 8\}\\
 &\geq\exp(-b_0(\ve,N)(\frac 1{c_0}-a+\la))= \ve^{1-c_0\hat a}\nonumber
 \end{eqnarray}
 where we choose $\la>0$ so small that $\hat a=a-\la>0$. By (\ref{3.21}) and
 (\ref{3.22}) we obtain
 \begin{equation}\label{3.23}
 \prod_{0\leq j<[T/N\ve][\ln^\ve]^{-1}}P\big\{\rho(\hat W^{\ve,N}_{0,j
 [\ln^2\ve]},\,\gam)\geq\del/4\big\}\leq(1-\ve^{1-c_0})^{T(N\ve\ln^2\ve)
 ^{-1}}.
 \end{equation}
 
 Taking $\ve_n=\frac 1n$ it follows from the Borel-Cantelli lemma together
 with the estimates (\ref{3.14}), (\ref{3.16}), (\ref{3.18}), (\ref{3.19})
 and (\ref{3.23}) that with probability one
 \begin{equation}\label{3.24}
 \limsup_{n\to\infty}\min_{0\leq l\leq[\frac Tn/N]}\rho(\hat W_t^{1/n,N},\gam)
 \leq\del.
 \end{equation}
 If $1/n\leq\ve\leq 1/(n-1)$ then $[tn]-[t/\ve]\leq T+1$ and
 \begin{equation}\label{3.25}
 |b_0(\ve,N,u)-b_0(\frac 1n,N,u)|\leq c_0\ln(\frac n{n-1})\to 0\quad\mbox{as}
 \quad n\to\infty.
 \end{equation}
 This together with (\ref{3.14}) and (\ref{3.24}) yields
 \begin{equation}\label{3.26}
 \limsup_{\ve\to 0}\inf_{0\leq t\leq T}\rho(W^{\ve,N}_t,\gam)\leq\del\quad
 \mbox{a.s.}
 \end{equation}
 Since $\Phi_0(1/c_0)$ is a compact set we can choose there a $\del$-net
 $\gam_1,\gam_2,...,\gam_{k(\del)}$ and then with probability one (\ref{3.26})
 will hold true simultaneously for all $\gam_i,\, i=1,...,k(\del)$ in place of
 $\gam$ there. It follows
 then that with probability one
 \begin{equation*}
 \limsup_{\ve\to 0}\sup_{\gam\in\Phi_0(1/c_0)}\inf_{0\leq t\leq T}
 \rho(W^{\ve,N}_t,\gam)\leq 2\del
 \end{equation*}
 and since $\del>0$ is arbitrary we obtain (\ref{2.11}) for $W^{\ve,N}_t$ in
 place of $V^{\ve,N}_t$ which, as explained at the beginning of this subsection
 gives (\ref{2.11}) and completes the proof of Theorem \ref{thm2.1}. \qed
 
 \subsection{Proof of Corollaries \ref{cor2.2} and \ref{cor2.3}}
 Under the conditions of Corollary \ref{cor2.2} there is no dependence on $t$
 of $S_t=S$ in (\ref{2.8}) and we consider
 \[
 \Phi(1/c)=\{\gam\in C([0,1],\bbR^d):\,\gam_0=0,\, S(\gam)\leq 1/c\}
 \]
 in (\ref{2.11}). Thus there is no dependence on $N$ of quantities in
 (\ref{2.10}) and (\ref{2.11}), so that the limit in $N$ is not relevant now.
 It follows from (\ref{2.10}) that all limit points as 
 $\ve\to 0$ of curves from $\cW^\ve_c$ belong to the compact set $\Phi(1/c)$.
 Now observe that (\ref{2.11}) means that with probability one any $\gam\in
 \Phi(1/c)$ is a limit point as $\ve\to 0$ of curves from $\cW^\ve_c$ which
 yields (\ref{2.12}). \qed
 
 In order to derive Corollary \ref{cor2.3} observe that (\ref{2.10}) implies,
 in particular, that for any continuous function $f$ on the space of
 curves $[0,1]\to\bbR^d$ with probability one,
 \begin{equation}\label{3.27}
 \lim_{N\to\infty}\limsup_{\ve\to 0}\sup_{0\leq t\leq T}f(V^{\ve,N}_t)\leq
 \limsup_{N\to\infty}\sup_{0\leq t\leq T}\sup_{\gam\in\Phi_{\tau(t,N)}
 (c^{-1}_{\tau(t,N)})}f(\gam).
 \end{equation}
 Now set $c_t=1/I_t(\be)$ assuming that $I_t(\be)<\infty$ for all $t\in[0,1]$.
 Since now $d=1$ we can define $f(\gam)=\gam(1),\,\gam(u)=\gam_u$. Then
 \begin{eqnarray}\label{3.28}
 &\limsup_{N\to\infty}\sup_{0\leq t\leq T}\sup_{\gam\in\Phi_{\tau(t,N)}
 (c^{-1}_{\tau(t,N)})}f(\gam)\\
 &=\limsup_{N\to\infty}\sup_{0\leq t\leq T}
 \{\gam(1):\,\gam\in\Phi_{\tau(t,N)}(I_{\tau(t,N)}(\be))\}=\be.\nonumber
 \end{eqnarray}
 Indeed, by convexity of each rate function $I_t$ for any 
 $\gam\in\Phi_{s}(I_{s}(\be))$,
 \[
 I_s(\be)\geq S_s(\gam)=\int_0^1I_s(\dot\gam(u))du\geq I_s(\int_0^1\dot\gam(u)
 du)=I_s(\gam(1))
 \]
 and by monotonicity of $I_s$ (see Appendix), $\be\geq\gam(1)$. On the other
 hand, take $\gam(u)=u\be,\, u\in[0,1],$ then $S_s(\gam)=I_s(\be)$ for all 
 $s\in[0,T]$ and $\gam(1)=\be$ implying (\ref{3.28}). Observe also that, in
 particular, if $\gam\in\Phi_0(I_0(\be))$ then by (\ref{2.11}) with probability
 one,
 \[
 \lim_{N\to\infty}\limsup_{\ve\to 0}\inf_{0\leq t\leq T}|V^\ve_{\tau(t,N)}(1)
 -\be|=0
 \]
 which together with (\ref{3.27}) and (\ref{3.28}) yields (\ref{2.13}).  \qed

\section{Continuous time case}\label{sec4}\setcounter{equation}{0}
\subsection{Basic estimates}\label{subsec4.1}
In view of (\ref{1.3}) we have 
\begin{eqnarray}\label{4.1}
&\ve^{-1}(X^\ve_{t/\ve+b_t(\ve,N,u)}-X^\ve_{t/\ve})=\int_{t/\ve}^
{t/\ve+b_t(\ve,N,u)}B(X^\ve_s,\xi_s)ds\\
&=\int_{t/\ve}^{t/\ve+b_t(\ve,N,u)}B_t(\xi_s)ds+\Psi^{(1)}_{\ve,N}
(t,u)+\Psi^{(2)}_{\ve,N}(t,u)\nonumber
\end{eqnarray}
where
\[
\Psi^{(1)}_{\ve,N}(t,u)=\int_{t/\ve}^{t/\ve+b_t(\ve,N,u)}
(B(X^\ve_s,\xi_s)-B(\bar X^\ve_s,\xi_s))ds
\]
and
\[
\Psi^{(2)}_{\ve,N}(t,u)=\int_{t/\ve}^{t/\ve+b_t(\ve,N,u)}
(B(\bar X^\ve_s,\xi_s)-B_t(\xi_s))ds.
\]
Similarly to (\ref{3.2}) and (\ref{3.3}) by (\ref{2.1}) and the averaging
principle (\ref{1.5}),
\begin{equation}\label{4.2}
|\Psi^{(1)}_{\ve,N}(t,u)|\leq L_1T\ln(1/\ve)\sup_{0\leq t\leq T(\ve^{-1}
-\ln\ve)}|X^\ve_t-\bar X^\ve_t|\to 0\,\,\,\mbox{as}\,\,\,\ve\to 0.
\end{equation}
Similarly to (\ref{3.4}) and (\ref{3.5}) by (\ref{1.4}), (\ref{1.6}) and
 (\ref{2.1}),
\begin{equation}\label{4.3}
|\Psi^{(2)}_{\ve,N}(t,u)|\leq L_1^2T^2\ve\ln^2(1/\ve)\to
 0\,\,\,\mbox{as}\,\,\,\ve\to 0.
\end{equation} 

It follows from (\ref{4.1})--(\ref{4.3}) that (\ref{2.10}) and (\ref{2.11})
will hold true in the continuous time case provided these limits are verified
for
\begin{equation}\label{4.4}
\hat V_t^{\ve,N}(u)=\frac 1{b_t(\ve,N)}\int_{t/\ve}^{t/\ve+
b_t(\ve,N,u)}(B_t(\xi_s)-\bar B_t)ds,\,\, t\in[0,T],\, u\in[0,1]
\end{equation}
in place of $V_t^{\ve,N}(u)$. Since in the continuous time case we have the
 same estimates as in (\ref{3.6}) it follows that, again, (\ref{2.10}) and
 (\ref{2.11}) will hold true here provided they are obtained for
\begin{equation}\label{4.5}
W_t^{\ve,N}(u)=\frac 1{b_t(\ve,N)}\int_{[t/\ve]}^{[t/\ve]+b_t(\ve,N,u)}
(B_{\tau(t,N)}(\xi_s)-\bar B_{\tau(t,N)})ds,
\end{equation}
$t\in[0,T],\, u\in[0,1]$ in place of $V_t^{\ve,N}(u)$.

\subsection{The upper bound}\label{subsec4.2}
Set
\begin{equation}\label{4.6}
\hat W^{\ve,N}_{k,l}(u)=b^{-1}_{kT/N}(\ve,N)\int_{l}^{l+b_{kT/N}(\ve,N)}
(B_{kT/N}(\xi_s)-\bar B_{kT/N})ds.
\end{equation}
In order to prove (\ref{2.10}) for $W_t^{\ve,N}$ in place of $V_t^{\ve,N}$
it suffices to obtain (\ref{3.8}) for $\hat W^{\ve,N}_{k,l}$ defined by 
(\ref{4.6}). The proof of this proceeds almost verbatim as for the discrete
time in Section \ref{subsec3.2} using now the continuous time case upper
large deviations bounds for normalized integrals (\ref{2.20}) taking into
account that $\hat W^{\ve,N}_{k,0}(u)=Y_{kTN^{-1},c_{kTN^{-1}}\ln(1/\ve)}(u)$.

\subsection{The lower bound}\label{subsec4.3}
Set
\[
\hat W_t^\ve(u)=\hat W_t^\ve(u,\om)=b^{-1}_0(\ve)\int_{[t/\ve]}^{[t/\ve]
+b_0(\ve,u)}(B_0(\xi_s(\om))-\bar B_0)ds,\, u\in[0,1]
\]
where $b_0(\ve,u)=b_0(\ve,N,u)$ and $b_0(\ve)=b_0(\ve,N)$ do not depend on
$N$. Clearly, for any $\gam\in C([0,1],\bbR^d)$,
\begin{equation}\label{4.7}
\inf_{0\leq t\leq T}\rho(W_t^{\ve,N},\gam)\leq\inf_{0\leq t<T/N}
\rho(W^{\ve,N}_t,\gam)=\inf_{0\leq t<T/N}\rho(\hat W^\ve_t,\gam).
\end{equation}
Introduce,
\[
U^\ve(k,\hat\om)=b^{-1}_0(\ve)\sum_{0\leq n\leq k}\hat B_0\circ\te^n(\hat\om)=
b^{-1}_0(\ve)\int_0^{\sig_k(\hat\om)}B_0(\xi_s(\hat\om))ds
\]
where $\hat B_0(\hat\om)=\hat B(\bar Z_0,\hat\om)$, $\sig_k=\sum_{i=0}^{k-1}
\vsig\circ\te^i$ and $\hat B$ is the same as in (\ref{2.17}). Set
\[
\up^\ve(u,\hat\om)=\max\{ j\geq 0:\,\sig_j(\hat\om)\leq b_0(\ve,u)\}.
\]
Observe that if $\sig_l(\hat\om)\leq[t/\ve]\leq\sig_{l+1}(\hat\om)$ then
\begin{equation}\label{4.8}
\sup_{0\leq u\leq 1}|\hat W_t^\ve(u,\hat\om)-U^\ve(\up^\ve(u,\te^l\hat\om),\te^l
\hat\om)|\leq 4b^{-1}_0(\ve)L_1L_3.
\end{equation}
It follows that for any $\del>8b^{-1}_0(\ve)L_1L_3$,
\begin{eqnarray}\label{4.9}
&\{\hat\om:\,\inf_{0\leq t\leq T/N}\rho(\hat W^\ve_t(\hat\om),\gam)>\del\}\\
&\subset\{\min_{0\leq l\leq\frac T{2\ve N\bar\vsig}}\rho(\hat U^\ve\circ\te^l,
\gam)>\del/2\}\cup\{\hat\om:\,\sig_{[\frac T{2\ve N\bar\vsig}]}\geq
\frac T{\ve N}\}\nonumber
\end{eqnarray}
where $\hat U^\ve(u)=\hat U^\ve(u,\hat\om)=U^\ve(\up^\ve(u,\hat\om),\hat\om)$.
It follows from (\ref{2.19}) that
\begin{equation}\label{4.10}
\hat P\{\sig_{[\frac T{2\ve N\bar\vsig]}}\geq\frac T{\ve N}\}\leq\ka^{-1}
\exp(-\frac \ka{\ve N})
\end{equation}
for some $\ka=\ka_T>0$.

Next, we deal with the other event in the right hand side of (\ref{4.9}).
Set
\[
\up_M^\ve(u,\hat\om)=\max\{ j\geq0:\,\sig_{j,M}(\hat\om)\leq b_0(\ve,u)\}
\]
where
\[
\sig_{j,M}=\sum_{i=0}^{j-1}\vsig_{i,M},\,\sig_{0,M}=0\,\,\mbox{and}\,\,
\vsig_{i,M}=E(\vsig\circ\te^i|\cF_{i-M,i+M}).
\]
Put also $\tilde B_{n,M}=\hat E(\hat B_0\circ\te^n|\cF_{n-M,n+M})$ and
define
\[
\hat U^\ve_{l,M}(u)=\hat U^\ve_{l,M}(u,\hat\om)=b^{-1}_0(\ve)\sum_
{l\leq n\leq l+\up^\ve_M(u,\hat\om)}\tilde B_{n,M}(\hat\om).
\]
Then
\begin{eqnarray}\label{4.11}
&|\hat U^\ve(u,\te^l\hat\om)-\hat U^\ve_{l,M}(u,\hat\om)|\leq L_1b^{-1}_0(\ve)
|\up^\ve(u,\hat\om)-\up^\ve_M(u,\hat\om)|\\
&+b^{-1}_0(\ve)\sum_{l\leq n\leq l+L_3b_0(\ve)}|\hat B_0\circ\te^n(\hat\om)-
\tilde B_{n,M}(\hat\om)|.\nonumber
\end{eqnarray}

Next, we estimate the right hand side of (\ref{4.11}). Observe that for all
$u\in[0,1]$,
\begin{equation}\label{4.12}
\{\hat\om:\, |\up^\ve_M(u,\hat\om)-\up^\ve(u,\hat\om)|\geq 2\}\subset
\cup_{k:\,\sig_k(\hat\om)\leq b_0(\ve)}\{\hat\om:\,|\sig_k(\hat\om)-\sig_{k,M}
(\hat\om)|>L^{-1}_3\}.
\end{equation}
Hence, by (\ref{2.17}) and the Chebyshev's inequality,
\begin{eqnarray}\label{4.13}
&\hat P\{\sup_{0\leq u\leq 1}|\up_M^\ve(u,\cdot)-\up^\ve(u,\cdot)|\geq 2\}
\leq\sum_{k:\,\sig_k(\hat\om)\leq b_0(\ve)}\hat P\{|\sig_k-\sig_{k,M}|\\
&>L^{-1}_3\}\leq L_3\sum_{k:\,\sig_k(\hat\om)\leq b_0(\ve)}\sum_{0\leq j\leq 
k-1}\hat E|\vsig_j-\vsig_{j,M}|\leq L^3_3b_0^2(\ve)\ka_3^{-1}e^{-\ka_3M}.
\nonumber\end{eqnarray}
By (\ref{2.16}) we also have that
\begin{equation}\label{4.14}
\sum_{l\leq n\leq l+L_3b_0(\ve)}\hat E|\hat B_0\circ\te^n-\tilde B_{n,m}|
\leq L_3b_0(\ve)\ka_3^{-1}e^{-\ka_3M},
\end{equation}
and so
\begin{eqnarray}\label{4.15}
&\hat P\{\rho(\hat U^\ve\circ\te^l,\hat U^\ve_{l,M})\geq\del/4\}\leq
\hat P\{\sup_{0\leq u\leq 1}|\up^\ve(u)-\up^\ve_M(u)|\\
&\geq L_1b_0(\ve)\del/8\}+\hat P\{\sum_{l\leq n\leq l+L_3b_0(\ve)}|
\hat B_0\circ\te^n-\tilde B_{n,M}|\geq b_0(\ve)\del/8\}\nonumber\\
&\leq (L^2_3b_0^2(\ve)+8/\del)\ka_3^{-1}e^{-\ka_3M}
\nonumber\end{eqnarray}
provided $b_0(\ve)\geq 16\del^{-1}L_1^{-1}$.

Observe that $\sig_{j,M}$ is $\cF_{-M,j+M}$-measurable and
\[
\{\up^\ve_M(u,\hat\om)=k\}=\{\sig_{k,M}\leq b_0(\ve,u)\}\cap\{\sig_{k+1,M}
>b_0(\ve,u)\}
\]
which is $\cF_{-M,k+M+1}$-measurable. Since always $\up^\ve_M(u,\hat\om)\leq 
L_3b_0(\ve)$ we obtain that $\hat U^\ve_{l,M}$ is 
$\cF_{l-M,l+M+[L_3b_0(\ve)]+1}$-measurable. Now we choose $M=M(\ve)=[\ln^2\ve]$
and obtain by (\ref{4.15}) that
\begin{eqnarray}\label{4.16}
&\hat P\{\min_{0\leq l\leq\frac T{2\ve n\bar\vsig}}\rho(\hat U^\ve\circ\te^l,
\gam)>\del/2\}\\
&\leq\hat P\{\min_{0\leq l\leq\frac T{6\ve N\bar\vsig M(\ve)}}
\rho(\hat U^\ve_{3lM(\ve),M(\ve)},\gam)>\del/4\}\nonumber\\
&+T(6\ve N\bar\vsig M(\ve))^{-1}L_3(L_3b^2_0(\ve)+8/\del)\ka_3^{-1}
e^{-\ka_3M(\ve)}.
\nonumber\end{eqnarray}

Introduce the event
\[
C^\ve_l=\{\rho(\hat U^\ve_{3lM(\ve),M(\ve)},\gam)>\del/4\}
\]
which is $\cF_{l-M(\ve),l+M(\ve)+[L_3b_0(\ve)]+1}$-measurable. Then
by (\ref{2.2}) for the probability $Q$ in the same way as in (\ref{3.17})
and (\ref{3.18}) it follows that,
\begin{eqnarray}\label{4.17}
&\hat P\big(\cap_{0\leq l\leq\frac T{6\ve N\bar\vsig M(\ve)}}C^\ve_l\big)\leq
L_2Q\big(\cap_{0\leq l\leq\frac T{6\ve N\bar\vsig M(\ve)}}C^\ve_l\big)\\
&\leq L_2\prod_{0\leq l\leq\frac T{6\ve N\bar\vsig M(\ve)}}Q(C_l^\ve)+
L_2\ve^{-1}\al([M(\ve)/2])\nonumber\\
&=L_2\prod_{0\leq l\leq\frac T{6\ve N\bar\vsig M(\ve)}}(1-Q(\hat\Om\setminus
C^\ve_l))+L_2\ve^{-1}\al([M(\ve)/2])\nonumber
\end{eqnarray}
provided $\ve>0$ is small enough.

Now,
\begin{equation}\label{4.18}
Q(\hat\Om\setminus C^\ve_l)\geq Q(\hat C^\ve_l)-L_2L_3(L^2_3b^2_0(\ve)+
16\del^{-1})\ka_3^{-1}e^{-\ka_3M(\ve)}
\end{equation}
where $\hat C^\ve_l=\{\rho(\hat U^\ve\circ\te^l,\gam)\leq\del/8\}$ and we use
(\ref{4.15}) with $\del/2$ in place of $\del$. Since $\te$ preserves the
measure $Q$,
\begin{equation}\label{4.19}
Q(\hat C^\ve_l)=Q(\hat C^\ve_0)\geq L^{-1}_2\hat P(\hat C_0^\ve).
\end{equation}
By (\ref{4.8}) for all $\hat\om$,
\begin{equation}\label{4.20}
\sup_{0\leq u\leq 1}|\hat U^\ve(u,\hat\om)-\hat W^\ve_0(u,\hat\om)|\leq
\frac {6L_1L_2}{b_0(\ve)}.
\end{equation}
Set $D^\ve=\{\hat\om:\,\rho(\hat W_0^\ve(\hat\om),\gam)\leq\del/16\}$ then
\begin{equation}\label{4.21}
\hat P(\hat C^\ve_0)\geq\hat P(D^\ve)\quad\mbox{provided}\quad \del>96L_1L_3
b_0^{-1}(\ve).
\end{equation}

Set $\hat D^\ve=\{\om:\,\rho(\hat W^\ve_0(\om),\gam)\leq\del/32\}$. Recal that
if $\hat\om=\vf\om$ then $\om=\vt^s\hat\om$ for some $0\leq s<\vsig(\hat\om)$,
and so $\xi_t(\om)=\xi_{t+s}(\hat\om)$ for any $t\geq 0$. This together with
(\ref{2.1}) and (\ref{2.18}) yields
\begin{equation}\label{4.22}
\sup_{0\leq u\leq 1}|\hat W^\ve_0(u,\om)-\hat W^\ve_0(u,\hat\om)|\leq 4L_1L_3
b_0^{-1}(\ve).
\end{equation}
Hence, if $\ve$ is small enough then
\begin{equation}\label{4.23}
\vf\hat D^\ve\subset D^\ve.
\end{equation}

In the same way as in the discrete time case we argue that since 
$\gam\in\Phi_0^{1/c_0}$ then the curve
\[
\eta_u=(1-\del(72\sup_{v\in[0,1]}|\gam_u|)^{-1})\gam_u,\,\, u\in[0,1]
\]
satisfies
\begin{equation}\label{4.24}
\rho(\gam,\eta)\leq\frac \del{72}\,\,\,\mbox{and}\,\,\, S_0(\eta)\leq S_0(\gam)
-a\leq\frac 1{c_0}-a
\end{equation}
for some $a>0$. This relies, again, on the strict monotonicity of the rate
function of large deviations in the domain where it is finite (see Appendix).
Next, we apply the lower large deviations bound in (\ref{2.9}) to
$Y_{0,c_0|\ln\ve|}(u)=\hat W_0^\ve(u)$ obtaining that for any $\la,\del>0$ 
there exists $\ve_0>0$ such that whenever $\ve\leq\ve_0$,
\begin{equation}\label{4.25}
\hat P(D^\ve)\geq P(\hat D^\ve)\geq\exp(-b_0(\ve)(\frac 1{c_0}-a+\la))=
\ve^{1-c_0\hat a}
\end{equation}
where we choose $\la>0$ so small that $\hat a=a-\la>0$.

Now, by (\ref{4.7}), (\ref{4.9}), (\ref{4.10}) and (\ref{4.22}) for $\ve$
small enough,
\begin{eqnarray}\label{4.26}
&P\{\inf_{0\leq t\leq T}\rho(W_t^{\ve,N},\gam)>2\del\}\leq 
P\{\inf_{0\leq t<T/N}\rho(\hat W_t^\ve,\gam)>2\del\}\\
&\leq\hat P\{\hat\om:\,\inf_{0\leq t<T/N}\rho(\hat W_t^\ve(\hat\om),\gam)>
\del\}\nonumber\\
&\leq\hat P\{\min_{0\leq l\leq T(2\ve N\bar\vsig)^{-1}}\rho(\hat U^\ve\circ
\te^l,\gam)>\del/2\}+\ka^{-1}\exp(-\ka/\ve N)\nonumber
\end{eqnarray}
for some $\ka>0$.
Next, by (\ref{4.16})--(\ref{4.19}), (\ref{4.21}), (\ref{4.25}) and 
(\ref{4.26}),
\begin{eqnarray}\label{4.27}
&P\{\inf_{0\leq t\leq T}\rho(W_t^{\ve,N},\gam)>2\del\}\\
&\leq L_2(1-\ve^{1-c_0\hat a}+L_4\del^{-1}e^{-\ka_4\ln^2\ve}\ln^2\ve)^{T(6
\ve N\bar\vsig\ln^2\ve)^{-1}}+L_4\ve^{-1}\del^{-1}e^{-\ka_4\ln^2\ve}
\nonumber\end{eqnarray}
for some $\ka_4,L_4>0$ independent of $\ve$.

Taking $\ve_n=1/n$ it follows from the Borel-Cantelli lemma that with
probability one,
\begin{equation}\label{4.28}
\limsup_{n\to\infty}\inf_{0\leq t\leq T}\rho(W_t^{1/n},\gam)\leq 2\del.
\end{equation}
If $1/n\leq\ve <\frac 1{n-1}$ then using (\ref{3.25}) we conclude again that
(\ref{4.28}) implies, in fact, that with probability one,
\begin{equation}\label{4.29}
\limsup_{\ve\to 0}\inf_{0\leq t\leq T}\rho(W_t^{\ve},\gam)\leq 2\del.
\end{equation}
Concluding in the same way as in the discrete time case of Section 
\ref{subsec3.3} by choosing a $\del$-net in $\Phi_0(1/c_0)$ and taking into 
account that
 $\del>0$ is arbitrary we obtain (\ref{2.11}) for $W_t^{\ve,N}$ in place
of $V_t^{\ve,N}$ which, as explained at the beginning of this section gives
(\ref{2.11}) and completes the proof of Theorem \ref{thm2.4}.   \qed

\subsection{Proof of Corollary \ref{cor2.5}}\label{subsec4.4}

The proof of Corollary \ref{cor2.5} proceeds, essentially, in the same way as 
the proofs of Corollaries \ref{cor2.2} and \ref{cor2.3} relying on properties of
large deviations rate functionals for the continuous time case (see Appendix).

\section{Appendix}\label{sec5}\setcounter{equation}{0}
\subsection{Applications}\label{subsec5.1}
The main applications in the discrete time case of Theorem \ref{thm2.1}
concern Markov chains and some classes of dynamical systems such as Axiom
A diffeomorphisms, expanding transformations and topologically mixing 
subshifts of finite type. We will restrict ourselves to several main setups
to which our results are applicable rather than trying to describe most
general situations. First, let $\xi_n,\, n\geq 0$ be a time homogeneous
Markov chain on a Polish state space $M$ whose transition probability
$P(x,\Gam)=P\{\xi_1\in\Gam|\xi_0=x\}$ satisfies
\begin{equation}\label{5.1}
\ka\nu(\Gam)\leq P(x,\Gam)\leq\ka^{-1}\nu(\Gam)
\end{equation}
for some $\ka>0$, a probability measure $\nu$ on $M$ and any Borel set 
$\Gam\subset M$. Then $\xi_n,\, n\geq 0$ is exponentially fast $\psi$-mixing
with respect to the family of $\sig$-algebras $\cF_{m,n}=\sig\{\xi_k,\,
m\leq k\leq n\}$ generated by the process (see, for instance, \cite{IL}).
The strong Doeblin type condition (\ref{5.1}) implies geometric ergodicity
\[
\| P(n,x,\cdot)-\mu\|\leq\be^{-1}e^{-\be n},\,\be>0
\]
where $\|\cdot\|$ is the variational norm, $P(n,x,\cdot)$ is the $n$-step 
transition probability and $\mu$ is the unique invariant measure of
$\{\xi_n,\, n\geq 0\}$ which makes it a stationary process.

In this situation the limit (\ref{2.6}) exists (see Lemma 4.3 in Ch.7 of
\cite{FW}) and $\exp(\Pi_t(\frb))$ turns out to be the principal eigenvalue
of the positive operator
\[
Qf(x)=E_xf(\xi_1)\exp\big((\frb,G_t(\xi_1))\big)
\]
where $E_x$ is the expectation provided $\xi_0=x$ (see \cite{Ki1}, \cite{Ki5} 
and references there). It is well known (see
\cite{Na}, \cite{IL}, \cite{GH} and references there) that $\Pi_t(\frb)$
is convex and differentiable in $\frb$. Furthermore, the Hessian matrix
$\nabla^2_\frb\Pi_t(\frb)|_{\frb=0}$ is positively definite if and only if
for each $\frb\in\bbR^d,\,\frb\ne 0$ the limiting variance
\begin{equation}\label{5.2}
\sig_\frb^2=\lim_{n\to\infty}n^{-1}E\big(\sum_{i=0}^n(\frb,G_t(\xi_i))\big)^2
\end{equation}
is positive. The latter holds true unless there exists a representation
$(\frb,G_t(\xi_n))=g(\xi_n)-g(\xi_{n-1}),\, n=1,2,...$ for some Borel function
$g$ (see \cite{IL}).

In the discrete time dynamical systems case we consider $\xi_n=\xi_n(x)
=f^nx,\, n\geq 0$ 
where $f:M\to M$ is a $C^2$ Axion A diffeomorphism on a hyperbolic set or
a topologically mixing subshift of finite type or a $C^2$ expanding 
transformation. Here $\xi_n,\, n\geq 0$ is considered as a stationary
process on the probability space $(M,\cF,\mu)$ where $\cF$ is the Borel
$\sig$-algebra and $\mu$ is a Gibbs measure constructed by a H\" older 
continuous function (see \cite{Bow}). Then the process $\xi_n$ is exponentially
fast $\psi$-mixing (see \cite{Bow}) with respect to the family of (finite)
$\sig$-algebras generated by cylinder sets in the symbolic setup of subshifts
of finite type or with respect to the corresponding $\sig$-algebras constructed
via Markov partitions in the Axiom A and expanding cases. Existence of the
limit (\ref{2.6}) and its form was proved in \cite{Ki2}. Here $\Pi_t(\frb)$
turns out to be the topological pressure for the function $(\frb, G_t)+\vf$
where $\vf$ is the potential of the corresponding Gibbs measure. The 
differentiability
properties of $\Pi_t(\frb)$ in $\frb$ are well known and, again, the Hessian
matrix $\nabla^2_\frb\Pi_t(\frb)|_{\frb=0}$ is positively definite if and only if
for each $\frb\in\bbR^d,\,\frb\ne 0$ the limiting variance (\ref{5.2}) is
positive where the expectation should be taken with respect to the 
corresponding Gibbs measure (see \cite{PP}, \cite{GH}, \cite{Ki1}--\cite{Ki3}
and references there). The latter holds true unless there exists a
coboundary representation $(\frb,G_t)=g\circ f-g$ for some bounded Borel
function $g$.

Next, we discuss the continuous time case. Here $\xi_t,\, t\geq 0$ can be a
nondegenerate random evolution on a compact manifold $M$, in particular,
a nondegenerate diffusion there. The existence and the form of the limit
(\ref{2.6}) in this case is shown in \cite{Ki5}. By discretizing time the
problem is reduced to the discrete time process $\xi_n,\, n\geq 0$ which
is exponentially fast $\psi$-mixing and in this case the continuous time
does not pose additional difficulties provided we consider the $\sig$-algebras
$\cF_{m,n}$ generated by $\xi_i,\, m\leq i\leq n$. We observe that this case 
fits in our general continuous time scheme taking the projection $\vf$ to be the
identity map, $\hat\Om=\Om$ and $\vsig(\hat\om)\equiv 1$. If the information
about the process comes only at some random times then we arrive at a more 
general setup of Theorem \ref{thm2.4} though its main motivation comes
from dynamical systems as described below.

We deal now with continuous time dynamical systems, namely, with a $C^2$
Axiom A
flows $f^t:M\to M$ on a hyperbolic set considered with a Gibbs measure built
by a H\" older continuous function. Using Markov partitions such a flow can
be represented by means of a suspension construction (see \cite{BR}) with
a transformation $\te:\,\hat M\to\hat M$ on the bases of elements of the
Markov partition and a roof function $\vsig$ so that $f^{\vsig(\hat x)}\hat x
=\te\hat x$ for each $\hat x\in\hat M$. Here $M$ is identified with the
space $\tilde M=\{(s,\hat x):\,\hat x\in\hat M,\, 0\leq s<\vsig(\hat x)\}$
and $f^t(s,\hat x)=(s+t,\hat x)$ for $s+t<\vsig(\hat x)$. Furthermore,
$\te^n,\, n\geq 0$ on $\hat M$ turns out to be an exponentially fast
$\psi$-mixing discrete time dynamical system preserving a Gibbs measure $Q$ 
constructed by a H\" older continuous function while $f^t$ preserves the
measure $P$ such that 
\[
\int g(s,\hat x)dP(s,\hat x)=(1/\bar\vsig)\int_{\hat M}\int_0^{\vsig(\hat 
x)}g(s,\hat x)dsdQ(\hat x)
\]
where $\bar\vsig=\int\vsig dQ$. The roof function $\vsig$ turns out to be
H\" older continuous and bounded away from zero and infinity. Large deviations
estimates for sums $\sum_{i=1}^n\vsig\circ\te^i$ follow from \cite{Ki1} and
existence of the limit in (\ref{2.6}) and its form follow from \cite{Ki2}
and \cite{Ki5}. Again, $\Pi_s(\frb)$ is the topological pressure of the flow
$f^t$ for the function $(\frb,G_s)+\vf$, with $\vf$ being the potential of the
corresponding Gibbs measure, and the differentiability properties
of $\Pi_s(\frb)$ in $\frb$ are well known (see, for instance, \cite{PP} and
\cite{Con}). Similarly to the discrete time case the Hessian matrix
is positively definite if and only if all limiting variances
\begin{equation}\label{5.3}
\sig_\frb^2=\lim_{T\to\infty}T^{-1}\int\big(\int_0^T(\frb,G_s\circ f^u)du\big)^2dP
\end{equation}
are positive when $\frb\ne 0$ which holds true unless there exists a coboundary
representation $(\frb,G_s)=g\circ f^t-g$ for some $t$ and a bounded Borel 
function $g$.

\subsection{Some properties of rate functions}\label{subsec5.2}
We collect here few properties of rate functions of large deviations
which are essentially well known but hard to find in major books on large
deviations. First, observe that if $\Pi(\frb)$, $\frb\in\bbR^d$ is a twice
differentiable function such that $\Pi(0)=0,\,\nabla_\frb\Pi(\frb)|_{\frb=0}
=0$ then $\Pi(\frb)=o(|\frb|)$, and so
\begin{equation}\label{5.4}
I(\be)=\sup_\frb((\frb,\be)-\Pi(\frb))>0
\end{equation}
unless $\be=0$. Indeed, by the above
\[
I(\be)\geq\del|\be|^2-\Pi(\del\be)>0
\]
if $\del>0$ is small enough. Curiously, positivity of the rate function is not
discussed in several books on large deviations without which upper large
deviations bounds do not make much sense.

Next, assume, in addition, that $\Pi$ is convex and has a positively definite
at zero Hessian matrix $\nabla_\frb^2\Pi(\frb)|_{\frb=0}$. Then $\Pi(\frb)\geq
0$ for all $\frb\in\bbR^d$ and for some $\del_1,\del_2>0$,
\begin{equation}\label{5.5}
\Pi(\frb)\geq\del_1|\frb|\quad\mbox{provided}\quad |\frb|>\del_2.
\end{equation}
It follows that if $|\be|<\del_1$ then $\frb_\be=\arg\sup((\frb,\be)-\Pi(\frb))$ 
satisfies $|\frb_\be|\leq\del_2$ and, in particular, $I(\be)<\infty$, i.e.
$I(\be)$ is finite in some neighborhood of $0$. 

Next, under the above conditions on $\Pi$ suppose that $I(\be)<\infty$
for some $\be\ne 0$. Then 
\begin{equation}\label{5.6}
I((1+\del)\be)>I(\be)\quad\mbox{for any}\quad\del>0.
\end{equation}
Indeed, for any $\ve>0$ there exists $\frb_{\be,\ve}$ such that
\[
(\frb_{\be,\ve},\be)-\Pi(\frb_{\be,\ve})\geq I(\be)-\ve.
\]
Since $\Pi(\frb_{\be,\ve})\geq 0$ we have
\[
I((1+\del)\be)\geq (1+\del)(\frb_{\be,\ve},\be)-\Pi(\frb_{\be,\ve})\geq
I(\be)+\del(I(\be)-\ve)-\ve>I(\be)
\]
provided $\ve<\del(1+\del)^{-1}I(\be)$ yielding (\ref{5.6}).

In the Erd\H os-R\' enyi law type results it is important to know where
a rate function $I(\be)$ is finite. This issue is hidden inside the functional
form of Theorems \ref{thm2.1} and \ref{thm2.4} but appears explicitly in
Corollary \ref{cor2.3} and in the classical form (\ref{1.1}). The discussion
on finiteness of rate functions is hard to find in books on large deviations
though without studying this issue lower bounds there do not have much sense.
We start with the rate functional $J(\nu)$ of the second level of large 
deviations for occupational measures
\begin{equation}\label{5.7}
\zeta_n=\frac 1n\sum_{k=0}^{n-1}\del_{\xi_k}\,\,\,\mbox{or}\,\,\, \zeta_t=
\frac 1t\int_0^t\del_{\xi_s}ds
\end{equation}
in the discrete or continuous time cases, respectively, where $\del_x$ denotes
the unit mass at $x$ (see \cite{Ki1}). Explicit formulas for $J(\nu)$ are
known when $\xi_k$ is a Markov chain whose transition probability satisfies
(\ref{5.1}) and when $\xi_k=f^kx$ with $f$ being an Axiom A diffeomorphism,
expanding transformation or subshift of finite type. In the former case
(see \cite{DV}),
\begin{equation}\label{5.8}
J(\nu)=-\inf_{u>0,\,\mbox{\small continuous}}\int\ln(\frac {Pu}u)d\nu
\end{equation}
and in the latter case (see \cite{Ki1}),
\begin{equation}\label{5.9}
\gathered J(\nu)=\left\{
\aligned &-\int\varphi d\mu -h_\nu(f)
\ \text{if}\, \nu \, \text{is}\, f\text{-invariant,}\\
&\infty\ \text{otherwise}
\endaligned\right. \endgathered
\end{equation}
where $h_\nu(f)$ is the Kolmogorov--Sinai entropy of $f$ with respect to
 $\nu$ and $\vf$ is the potential of the corresponding Gibbs measure $\mu$
 playing the role of probability here. 
 
 In the continuous time case these functionals have explicit forms for
 diffusions $\xi_t$ (see \cite{DV}),
 \begin{equation}\label{5.10}
 J(\nu)=-\inf_{u>0,\,\mbox{\small is}\, C^2}\int\frac {Lu}ud\nu,
 \end{equation}
 where $L$ is the corresponding generator, and for Axiom A flows
 $\xi_t=f^tx$ where (see \cite{Ki1}),
 \begin{equation}\label{5.11}
\gathered J(\nu)=\left\{
\aligned &-\int\varphi d\mu -h_\nu(f^1)
\ \text{if}\, \nu \, \text{is}\, f\text{-invariant,}\\
&\infty\ \text{otherwise}
\endaligned\right. \endgathered
\end{equation}
with the same notations as in (\ref{5.9}).

Necessary and sufficient conditions for finiteness of $J(\nu)$ in the Markov
chain and diffusion cases are given in \cite{DV} while in the above
dynamical systems cases $J(\nu)<\infty$ for any $f$-invariant measure $\nu$.
If 
\begin{equation}\label{5.12}
\Pi(\frb)=\lim_{n\to\infty}\frac 1n\ln E\exp\big(\sum_{j=0}^{n-1}(\frb,
G(\xi_j))\big)
\end{equation}
in the discrete time case or
\begin{equation}\label{5.13}
\Pi(\frb)=\lim_{t\to\infty}\frac 1t\ln E\exp\big(\int_{0}^{t}(\frb,
G(\xi_s))ds\big)
\end{equation}
in the continuous time case, where $\xi_t$ is a stationary process as above
on a compact space $M$ and $G\not\equiv 0$ is a continuous vector function
with $EG(\xi_0)=0$, then by the contraction principle (see, for instance,
\cite{DZ}) the rate function $I(\be)$ given by (\ref{5.4}) can be represented
 as
 \begin{equation}\label{5.14}
 I(\be)=\inf\{ J(\nu):\,\int Gd\nu=\be\}
 \end{equation}
 where the infinum is taken over the space $\cP(M)$ of probability measures
 on $M$.
 
 Set 
 \[
 \Gam=\{\be\in\bbR^d:\,\exists\nu\in\cP(M)\,\,\mbox{such that}\,\,\int Gd\nu=
 \be \,\,\mbox{and}\,\, J(\nu)<\infty\}
 \]
 and let $Co(\Gam)$ be the interior of the convex hull of $\Gam$. Then
 \begin{equation}\label{5.15}
 I(\be)<\infty\,\,\mbox{for any}\,\,\be\in Co(\Gam).
 \end{equation}
 Indeed, any $\be\in Co(\Gam)$ can be represented as $\be=p_1\be_1+p_2\be_2$
 with $\be_1,\be_2\in\Gam$, $p_1,p_2\geq 0$ and $p_1+p_2=1$. Then $\be_1=
 \int Gd\nu_1,\,\be_2=\int Gd\nu$, and so $\int Gd\nu=\be$ for $\nu=p_1\nu_1
 +p_2\nu_2$. Since $J(\nu_1),\, J(\nu_2)<\infty$ then by convexity of $J$ we
 have that $J(\nu)\leq p_1J(\nu_1)+p_2J(\nu_2)<\infty$, and so (\ref{5.15})
 holds true.
 
 When $d=1$, i.e. when $G$ is (not vector) function we can give another
 description of the domain where $I(\be)<\infty$. In this case set
 \begin{equation}\label{5.16}
 \be_+=\sup\{\be:\,\be\in\Gam\}\,\,\mbox{and}\,\,\be_-=\inf\{\be:\be\in\Gam\}.
 \end{equation}
 Then by (\ref{5.15}), $I(\be)<\infty$ for any $\be\in (\be_-,\be_+)$. It is
 possible to extract from \cite{DK} that under $\psi$-mixing,
 \begin{equation}\label{5.17}
 \be_+=\lim_{n\to\infty}\frac 1ness\sup\sum_{j=0}^{n-1}G(\xi_j)\,\,\mbox{and}
 \,\, \be_-=\lim_{n\to\infty}\frac 1ness\inf\sum_{j=0}^{n-1}G(\xi_j).
 \end{equation}
 Since Axiom A flows are not $\psi$-mixing, in general, we will give another
 proof for this case.
 
 Let
 \begin{eqnarray}\label{5.18}
 &\be^*_+=\lim_{t\to\infty}\frac 1t\sup_x\int_0^tG\circ f^sds\,\,\mbox{and}\\
 &\be^*_-=\lim_{t\to\infty}\frac 1t\inf_x\int_0^tG\circ f^sds=-
 \lim_{t\to\infty}\frac 1t\sup_x(-\int_0^tG\circ f^sds).\nonumber
 \end{eqnarray}
 The limits in (\ref{5.17}) exist since $a(t)=\sup_x\int_0^tG\circ f^sds$ is
 subadditive $a(t+s)\leq a(t)+a(s)$. Since $G$ is a continuous function on a
 compact space $M$ we can find $x_t$ such that $a(t)=\int_0^tG\circ f^s(x_t)ds$.
 Consider the family of occupational measures
 \[
 \nu_t=\frac 1t\int_0^t\del_{f^sx_t}ds.
 \]
 Then any weak limit $\tilde\nu$ of $\nu_t$ as $t\to\infty$ is an 
 $f^t$-invariant
 measure and $\int gd\tilde\nu=\be^*_+$. It follows that $\be^*_+\leq\be_+$
 where $\be_+$ is given by (\ref{5.16}). On the other hand, $a(t)\geq t^{-1}
 \int_0^tG\circ f^s(x)ds$, and so for any $x$,
 \[
 \int Gd\tilde\nu\geq\lim\sup_{t\to\infty}\frac 1t\int_0^tG\circ f^s(x)ds.
 \]
 Hence, $\int Gd\tilde\nu\geq\int Gd\nu$ for any $f^t$-invariant probability
 measure $\nu$. Hence, $\be_+^*=\be_+$ and similarly we obtain that $\be_-^*=
 \be_-$.
 
 Even in the classical i.i.d. case of the Cram\' er theorem which is relevant 
 to the original form (\ref{1.1}) of the Erd\H os-R\' enyi law finiteness
 of the rate function is rarely discussed in details. Here, we provide a simple
 argumwnt. Let $\xi_1,\xi_2,...$ be i.i.d. random variables such that
 $E\xi_1=0$ and $\Pi(\frb)=\ln Ee^{\frb\xi_1}<\infty$ for all real $\frb$. Set
 \[
 I(\be)=\sup_\frb(\frb\be-\Pi(\frb)),\, \be_+=\|\xi_1^+\|_\infty=ess\sup\xi_1
 \,\,\mbox{and}\,\, \be_-=-\|\xi_1^-\|_\infty=ess\inf\xi_1.
 \]
  Then
 \begin{equation}\label{5.19}
 I(\be)<\infty\,\,\mbox{for any}\,\,\be\in(\be_-,\be_+)\,\,\mbox{and}\,\,
 I(\be)=\infty\,\,\mbox{if}\,\,\be\not\in[\be_-,\be_+].
 \end{equation}
 Indeed, if $0\leq\be<\be_+$ then $P\{\xi_1>\be\}=p_\be>0$. Hence,
 \begin{equation}\label{5.20}
 I(\be)=-\inf_{\frb\geq 0}\ln(e^{-\frb\be}Ee^{\frb\xi_1})\leq -\ln p_\be<\infty
 \end{equation}
 and similarly for $0\geq\be>\be_-$. On the other hand, if $\be>\be_+$ then
 \begin{equation}\label{5.21}
 I(\be)=-\inf_{\frb\geq 0}\ln(e^{-\frb\be}Ee^{\frb\xi_1})\geq -\inf_{\frb\geq 0}
 \frb(\be_+-\be)=\infty
 \end{equation}
 and similarly for $\be<\be_-$.

\bibliography{matz_nonarticles,matz_articles}
\bibliographystyle{alpha}

\end{document}